\documentclass{article}

\title{\LARGE \textbf{On $k$-ended spanning and dominating trees}}
\author{Zh.G. Nikoghosyan\footnote{G.G. Nicoghossian (up to 1997)}  }

\begin{document}

\maketitle

\begin{abstract}

A tree with at most $k$ leaves is called a $k$-ended tree. A spanning 2-ended tree is a Hamilton path. A Hamilton cycle can be considered as a spanning 1-ended tree. The earliest result concerning spanning trees with few leaves states that if $k$ is a positive integer and $G$ is a connected graph of order $n$ with $d(x)+d(y)\ge n-k+1$ for each pair of nonadjacent vertices $x,y$, then $G$ has a spanning $k$-ended tree. In this paper, we improve this result in two ways, and an analogous result is proved for dominating $k$-ended trees based on the generalized parameter $t_k$ - the order of a largest $k$-ended tree. In particular, $t_1$ is the circumference (the length of a longest cycle), and $t_2$ is the order of a longest path.  \\


\end{abstract}

\section{Introduction}

Throughout this article we consider only finite undirected graphs without loops or multiple edges. The set of vertices of a graph $G$ is denoted by $V(G)$ and the set of edges by $E(G)$.  A good reference for any undefined terms is $\cite{[1]}$.

For a graph $G$, we use $n$, $\delta$ and $\alpha$  to denote the order (the number of vertices), the minimum degree and the independence number of $G$, respectively. For a subset $S\subseteq V(G)$, we denote by $G[S]$ the subgraph of $G$ induced by $S$. If $\alpha\ge k$ for some integer $k$, let $\sigma_k$ be the minimum degree sum of an independent set of $k$ vertices; otherwise we let $\sigma_k=+ \infty$.

If $Q$ is a path or a cycle in a graph $G$, then the order of $Q$, denoted by $|Q|$, is $|V(Q)|$. Each vertex and edge in $G$ can be interpreted as simple cycles of orders 1 and 2, respectively. The graph $G$ is hamiltonian if $G$ contains a Hamilton cycle, i.e. a cycle containing every vertex of $G$. A cycle $C$ of $G$ is said to be dominating if $V(G-C)$ is an independent set of vertices.

We write a cycle $Q$ with a given orientation by $\overrightarrow{Q}$. For $x,y\in V(Q)$, we denote by $x\overrightarrow{Q}y$ the subpath of $Q$ in the chosen direction from $x$ to $y$. For $x\in V(Q)$, we denote the successor and the predecessor of $x$ on $\overrightarrow{Q}$ by $x^+$ and $x^-$, respectively.

A vertex of degree one is called an end-vertex, and an end-vertex of a tree is usually called a leaf. The set of end-vertices of $G$ is denoted by $End(G)$ . For a positive integer $k$, a tree $T$ is said to be a $k$-ended tree if $|End(T)|\le k$. A Hamilton path is a spanning 2-ended tree. A Hamilton cycle can be interpreted as a spanning 1-ended tree. In particular, $K_2$ is hamiltonian and is a 1-ended tree. We denote by $t_k$ the order of a largest $k$-ended tree in $G$. In particular, $t_1$ is the order of a longest cycle, and $t_2$ is the order of a longest path in $G$.

For two vertices $u$ and $v$ of $G$, let $d_G(u,v)$ denote the distance between $u$ and $v$. For a vertex $v$ of $G$, the distance between $v$ and a subgraph $H$ is defined to be the minimum value of $d_G(v,x)$ for all $x\in V(H)$, and denoted by $d_G(v,H)$. If $X\subseteq V(G)$ then $d_G(v,X)$ is defined analogously.

Our starting point is the earliest sufficient condition for a graph to be hamiltonian due to Dirac \cite{[3]}. \\

\noindent\textbf{Theorem A} \cite{[3]}. Every graph with $\delta\ge \frac{n}{2}$ is hamiltonian.  \\

In 1960, Ore \cite{[10]} improved Theorem A by replacing the minimum degree $\delta$ with the arithmetic mean $\frac{1}{2}\sigma_2$ of two smallest degrees among pairwise nonadjacent vertices.\\

\noindent\textbf{Theorem B} \cite{[10]}. Every graph with $\sigma_2\ge n$ is hamiltonian.\\

The analog of Theorem B for Hamilton paths follows easily.\\

\noindent\textbf{Theorem C} \cite{[10]}. Every graph with $\sigma_2\ge n-1$ has a Hamilton path.\\

In 1971, Las Vergnas \cite{[4]} gave a degree condition that guarantees that any forest in $G$ of limited size and with a limited number of leaves can be extended to a spanning tree of $G$ with a limited number of leaves in an appropriate sense. This result implies as a corollary a degree sum condition for the existence of a tree with at most $k$ leaves including Theorem B and Theorem C as special cases for $k=1$ and $k=2$, respectively.   \\

\noindent\textbf{Theorem D} \cite{[2]}, \cite{[4]}, \cite{[7]}. If $G$ is a connected graph with $\sigma_2\ge n-k+1$ and $k$ a positive integer, then $G$ has a spanning $k$-ended tree.\\

However, Theorem D was first openly formulated and proved in 1976 by the author \cite{[7]} and was reproved in 1998 by  Broersma and Tuinstra \cite{[2]}. Moreover,  the full characterization of connected graphs without spanning $k$-ended trees was given in \cite{[6]} when $\sigma_2\ge n-k$ including well-known characterization of connected graphs without Hamilton cycles when $\sigma_2\ge n-1$. This particular result was reproved in 1980 by Nara Chie \cite{[5]}.

In this paper we prove that the connectivity condition in Theorem D can be removed, and the conclusion can be strengthened.  \\

\noindent\textbf{Theorem 1}. If $G$ is a graph with $\sigma_2\ge n-k+1$ and $k$ a positive integer, then $G$ has a spanning $k$-ended forest.\\

The next improvement of Theorem D is based on parameter $t_k$ including the circumference (the length of a longest cycle) and the length of a longest path in a graph for $k=1$ and $k=2$, respectively.\\

\noindent\textbf{Theorem 2}. Let $G$ be a connected graph with $\sigma_2\ge t_{k+1}-k+1$ and $k$ a positive integer. Then $G$ has a spanning $k$-ended tree.\\

The graph $(\delta+k)K_1+K_\delta$ shows that the bound $t_{k+1}-k+1$ in Theorem 2 cannot be relaxed to $t_k-k+1$.

Finally, we give a degree sum condition for dominating $k$-ended trees.\\

\noindent\textbf{Theorem 3}. If $G$ is a connected graph with $\sigma_3\ge t_{k+1}-2k+4$ for some integer $k\ge2$, then $G$ has a dominating $k$-ended tree.\\

The graph $(\delta+k-1)K_2+K_{\delta-1}$ shows that the bound $t_{k+1}-2k+4$ in Theorem 3 cannot be relaxed to $t_k-2k+4$.

The following corollary follows immediately.\\

\noindent\textbf{Corollary 1}. If $G$ is a connected graph with $\sigma_3\ge n-2k+4$ for some integer $k\ge2$, then $G$ has a dominating $k$-ended tree.\\

The graph $(\delta+k-1)K_2+K_{\delta-1}$ shows that the bound $\sigma_3\ge t_{k+1}-2k+4$ in Theorem 3 cannot be relaxed to $\sigma_3\ge t_{k+1}-2k+3$.

We present also some earlier results concerning spanning $k$-ended trees that are not included in the recent survey paper \cite{[11]}. We call a graph $G$ hypo-$k$-ended if $G$ has no a spanning $k$-ended tree, but for any $v\in V(G)$, $G-v$ has a spanning $k$-ended tree.  \\

\noindent\textbf{Theorem E} \cite{[8]}. For each $k\ge 3$, the minimum number of vertices (edges, faces, respectively) of a simple 3-polytope without a spanning $k$-ended tree is $8+3k$ ($12+6k$, $6+3k$, respectively).\\

\noindent\textbf{Theorem F} \cite{[9]}. For each $n\ge 17k$ and $k\ge 2$, except possible for $n=17k+1$, $17k+2$, $17k+4$ and $17k+7$, there exist hypo-$k$-ended graphs of order $n$.

\section{Proofs}

\noindent\textbf{Proof of Theorem 1}. Let $G$ be a graph with
$\sigma_2\ge n-k+1$ and let $H_1,...,H_m$ be the connected components
of $G$. Let $\overrightarrow{P}=x\overrightarrow{P}y$ be a longest path in $H_1$.
If $|P|\ge n-k+2$ then $|G-P|=n-|P|\le k-2$,
 implying that $G$ has a spanning $k$-ended forest. Now let $|P|\le n-k+1$. Since $P$ is extreme, we have
 $N(x)\cup N(y)\subseteq V(P)$.  Recalling also that $\sigma_2\ge n-k+1$, we have (by standard arguments) $N(x)\cap N^+(y)\not=\emptyset$, implying that $G[V(P)]$ is hamiltonian. Further, if $|V(P)|<|V(H_1)|$
 then we can form a path longer than $P$, contradicting the maximality of $P$. Hence, $|V(P)|=|V(H_1)|$, that is $H_1$ is hamiltonian as well. By a similar argument, $H_i$ is hamiltonian for each $i\in \{1,...,m\}$ and therefore, has a spanning tree with exactly one leaf. Thus, $G$ has a spanning forest with exactly $m$ leaves.

It remains to show that $m\le k$. If $m=1$ then $G$ has a spanning 1-ended tree and therefore, has a spanning $k$-ended tree. Let $m\ge 2$ and let $x_i\in V(H_i)$ $(i=1,...,m)$. Clearly, $\{x_1,x_2,...,x_m\}$ is an independent set of vertices. Since $d(x_i)\le |V(H_i)|-1$, we have
$$
\sigma_2\le \sigma_m\le\sum_{i=1}^md(x_i)\le \sum_{i=1}^m|V(H_i)|-m=n-m.
$$
On the other hand, by the hypothesis, $\sigma_2\ge n-k+1$, implying that $m\le k-1$.  \qquad \rule{7pt}{6pt}\\

\noindent\textbf{Proof of Theorem 2}. Let $G$ be a connected graph with $\sigma_2\ge t_{k+1}-k+1$ for some positive integer $k$. \\

\textbf{Case 1}. $G$ is hamiltonian.

By the definition, $G$ has a spanning 1-ended tree $T_1$. Since $k\ge 1$, $T_1$ is a spanning $k$-ended tree.\\

\textbf{Case 2}. $G$ is not hamiltonian.

Let $T_2$ be a longest path in $G$.\\

\textbf{Case 2.1}. $\sigma_2\ge t_2$.

By standard arguments, $G[V(T_2)]$ is hamiltonian. If $t_2<n$ then recalling that $G$ is connected, we can form a path longer than $T_2$, contradicting the maximality of $T_2$. Otherwise $G$ is hamiltonian and we can argue as in Case 1.\\

\textbf{Case 2.2}. $\sigma_2\le t_2-1$.

If $k=1$ then by the hypothesis, $\sigma_2\ge t_2$, implying that $G$ is hamiltonian and we can argue as in Case 1. Let $k\ge 2$. Extend $T_2$ to a $k$-ended tree $T_k$ and assume that $T_k$ is as large as possible. If $T_k$ is a spanning tree then we are done. Let $T_k$ is not spanning. Then $|End(T_k)|=k$ since otherwise we can form a new $k$-ended tree larger than $T_k$, contradicting the maximality of $T_k$. Now extend $T_k$ to a largest $(k+1)$-ended tree $T_{k+1}$. Recalling that $T_k$ is a largest $k$-ended tree, we get $|End(T_{k+1})|=k+1$ and therefore,
$$
t_{k+1}\ge |T_{k+1}|= |T_2|+|T_{k+1}-T_2|.
$$
Observing that $|T_2|=t_2$ and $|T_{k+1}-T_2|\ge |End(T_{k+1})|-2=k-1$, we get 
$$
t_{k+1}\ge t_2+k-1\ge \sigma_2+k,
$$
contradicting the hypothesis.        \qquad \rule{7pt}{6pt}\\

\noindent\textbf{Proof of Theorem 3}. Let $G$ be a connected graph with $\sigma_3\ge t_{k+1}-2k+4$ for some integer $k\ge 2$, and let $\overrightarrow{T_2}=x\overrightarrow{T_2}y$ be a longest path in $G$. If $T_2$ is a dominating path then we are done. Otherwise, since $G$ is connected,  we can choose a path  $\overrightarrow{Q}=w\overrightarrow{Q}z$ such that $V(T_2\cap Q)=\{w\}$ and $|Q|\ge 3$. Assume that $|Q|$ is as large as possible. Put $T_3=T_2\cup Q$. Since $T_2$ and $Q$ are extreme, we have $N(x)\cup N(y)\subseteq V(T_2)$ and $N(z)\subseteq V(T_3)$. Let $w^+$ be the successor of $w$ on $T_2$. If $xy\in E$ then $T_3+xy-w^+w$ is a path longer than $T_2$, a contradiction. Let $xy\not\in E$. By the same reason, we have $xz,yz\not\in E$. Thus, $\{x,y,z\}$ is an independent set of vertices.\\

\noindent\textbf{Claim 1}. $N^-(x)\cap N^+(y)\cap N(z)=\emptyset$.

\noindent\textbf{Proof}. Assume the contrary.\\
 
\textbf{Case 1}. $v\in N^-(x)\cap N^+(y)$. 

If $v=w$ then $xv^+\in E$ and $T_3+xv^+-vv^+$ is a path longer than $T_2$, a contradiction. Suppose without loss of generality that $v\in V(w^+\overrightarrow{T_2}y)$. If $v=w^+$ then $T_3+xv^+-wv-vv^+$ is a path longer than $T_2$, a contradiction. Finally, if $v\in V(w^{+2}\overrightarrow{T_2}y)$ then
$$
T_3+xv^++yv^--vv^--vv^+-ww^+
$$
is a path longer than $T_2$, a contradiction.\\

\textbf{Case 2}. $v\in N^-(x)\cap N(z)$. 

If $v\in V(x\overrightarrow{T_2}w^{-2})$ then
$$
T_2+xv^++zv-vv^+-ww^-
$$
is a path longer than $T_2$, a contradiction. Next, if $v=w^-$ then $T_2+zw^--ww^-$ is a path longer than $T_2$, a contradiction. Further, if $v=w$ then $T_2+xv^+-ww^+$ is a path longer than $T_2$, a contradiction. Finally, if $v\in V(w^+\overrightarrow{T_2}y)$ then
$$
T_2+xv^++zv-ww^+-vv^+
$$
is a path longer than $T_2$, a contradiction.\\

\textbf{Case 3}. $v\in N^+(y)\cap N(z)$.

By a symmetric argument, we can argue as in Case 2. Claim 1 is proved.   \ \ \ \   $\triangle$ \\

By Claim 1,
$$
t_3\ge |T_3|\ge |N^-(x)|+|N^+(y)|+|N(z)|+|\{z\}|
$$
$$
=d(x)+d(y)+d(z)+1\ge \sigma_3+1.  \eqno{(1)}
$$
If $k=2$ then by the hypothesis, $\sigma_3\ge t_{k+1}-2k+4=t_3$, contradicting (1). Let $k\ge 3$. If $T_3$ is a dominating 3-ended tree then clearly we are done. Otherwise $G-T_3$ contains an edge and we can extend $T_3$ to a largest 4-ended tree $T_4$ with $|T_4|\ge |T_3|+2$. If $k=3$, then by the hypothesis, $\sigma_3\ge t_{k+1}-2k+4=t_4-2$. On the other hand, by (1), $t_4\ge |T_4|\ge |T_3|+2\ge \sigma_3+3$, a contradiction. Hence $k\ge 4$. If $T_4$ is dominating, then we are done. Otherwise we can extend $T_4$ to a largest 5-ended tree $T_5$ with $|T_5|\ge |T_4|+2\ge |T_3|+4$. This procedure may be repeated until a dominating $(r+1)$-ended tree $T_{r+1}$ is found. If $r+1\le k$ then we are done. Let $r\ge k$. Then
$$
t_{k+1}\ge |T_{k+1}|\ge |T_3|+2(k-2)
$$
$$
\ge \sigma_3+2k-3\ge t_{k+1}+1,
$$
a contradiction. The proof is complete.  \qquad \rule{7pt}{6pt}\\

\noindent Institute for Informatics and Automation Problems\\ National Academy of Sciences\\
P. Sevak 1, Yerevan 0014, Armenia\\
E-mail: zhora@ipia.sci.am  

\end{document}